\documentclass[conference, a4paper]{IEEEtran}
\IEEEoverridecommandlockouts

\ifCLASSINFOpdf
\usepackage[pdftex]{graphicx}
\DeclareGraphicsExtensions{.pdf,.jpeg,.png}
\else
\fi

\usepackage{cite}
\usepackage{amsmath,amssymb,amsfonts}
\usepackage{algorithm,algorithmic}
\usepackage{graphicx}
\usepackage{textcomp}
\usepackage{xcolor}
\usepackage[left=1.4cm, right=1.4cm, top=1.7cm]{geometry}
\usepackage{anyfontsize}
\usepackage{csquotes}
\usepackage{amsthm}
\usepackage{mathtools}
\usepackage{bm}
\usepackage{enumitem}
\usepackage{flushend}
\usepackage{bbm}
\usepackage{hyperref}
\usepackage{multirow}
\usepackage{cancel}
\usepackage{array}
\usepackage{multicol}

\makeatletter
\let\oldlt\longtable
\let\endoldlt\endlongtable
\def\longtable{\@ifnextchar[\longtable@i \longtable@ii}
\def\longtable@i[#1]{\begin{figure}[t]
		\onecolumn
		\begin{minipage}{0.5\textwidth}
			\oldlt[#1]
		}
		\def\longtable@ii{\begin{figure}[t]
				\onecolumn
				\begin{minipage}{0.5\textwidth}
					\oldlt
				}
				\def\endlongtable{\endoldlt
				\end{minipage}
				\twocolumn
		\end{figure}}
		\makeatother

		\DeclareMathOperator{\argmin}{argmin}
		
		\DeclareSymbolFont{ugmL}{OMX}{mdugm}{m}{n}
		\SetSymbolFont{ugmL}{bold}{OMX}{mdugm}{b}{n}
		\DeclareMathAccent{\wideparen}{\mathord}{ugmL}{"F3}

		\newtheorem{defn}{Definition}

		\usepackage{booktabs}
		\usepackage{xcolor}
		\usepackage{subcaption}
		
		\makeatletter
		\newcommand*\titleheader[1]{\gdef\@titleheader{#1}}
		\AtBeginDocument{%
			\let\st@red@title\@title
			\def\@title{%
				\bgroup\normalfont\normalsize\centering\@titleheader\par\egroup
				\vskip0.2em\st@red@title}
		}
		\makeatother
		
		\makeatletter
		\renewcommand{\fnum@figure}{Figure \thefigure}
		\makeatother

\title{{Optimization-Guided Exploration of Advanced Air Mobility Congestion Management Strategies with Stochastic Demands} \\
\thanks{\noindent
\textsuperscript{$\dagger$}corresponding author: \href{mailto:maxzli@umich.edu}{maxzli@umich.edu}.\\This research project was supported by MITRE's Independent Research \& Development Program.}
}

\titleheader{First US-Europe Air Transportation Research and Development Symposium (ATRDS2025)}

\author{\IEEEauthorblockN{Haochen Wu\textsuperscript{1}, Lesley A. Weitz\textsuperscript{2}, Jeffrey M. Henderson\textsuperscript{2}, Max Z. Li\textsuperscript{1,$\dagger$}}
\IEEEauthorblockA{\textsuperscript{1}University of Michigan, Ann Arbor, MI, USA\\
\textsuperscript{2}The MITRE Corporation, McLean, VA, USA\\
\{\href{mailto:haocwu@umich.edu}{haocwu}, \href{mailto:maxzli@umich.edu}{maxzli}\}@umich.edu,
\{\href{mailto:lweitz@mitre.org}{lweitz}, \href{mailto:jmhenderson@mitre.org}{jmhenderson}\}@mitre.org}
}

\renewcommand\IEEEkeywordsname{Keywords}
\IEEEaftertitletext{\vspace{-1\baselineskip}}

\begin{document}

\maketitle


\noindent \begin{abstract}
Advanced Air Mobility (AAM) represents an evolution of the air transportation system by introducing low-altitude, potentially high-traffic environments. AAM operations will be enabled by both new aircraft, as well as new safety- and efficiency-critical supporting infrastructure. Published concepts of operations from both public and private sector entities establish notions such as federated management of the airspace and public-private partnerships for AAM air traffic, but there is a gap in the literature in terms of integrated tools that consider all three critical elements: AAM fleet operators (\emph{lower} layer), airspace service providers (\emph{middle} layer), and overall system governance from the legacy air navigation service provider (\emph{upper} layer). In this work, we explore modeling congestion management within the AAM setting using a bi-level optimization approach, focusing on (1) time-varying, stochastic AAM demand, (2) differing congestion management strategies, and (3) the impact of unscheduled, \enquote{pop-up} demand. We show that our bi-level formulation can be tractably solved using a Neural Network-based surrogate which returns solution qualities within 0.1-5.2\% of the optimal solution. Additionally, we show that our congestion management strategies can reduce congestion by 25.7-39.8\% when compared to the scenario of no strategies being applied. Finally, we also show that while pop-up demand degrades congestion conditions, our congestion management strategies fare better against pop-up demand than the no strategy scenario. The work herein contributes a rigorous modeling and simulation tool to help evaluate future AAM traffic management concepts and strategies.
\end{abstract}



\section{Introduction}	\label{sec:intro}

Advanced Air Mobility (AAM) is a future vision for commonplace air transport in lower altitude, higher operational density settings using new aviation technologies. The U.S. AAM Coordination and Leadership Act \cite{USAAMCLA} provides one of several formal definitions, defining AAM as \enquote{a transportation system that moves people and property by air between two points in the U.S.~using aircraft with advanced technologies, including electric aircraft, or electric vertical takeoff and landing (eVTOL) aircraft, in both controlled and uncontrolled airspace.} Other definitions specify that AAM should provide transportation for people and cargo where it is not currently or easily served by surface transportation or existing aviation modes \cite{Texas_AAM_2024}. AAM has the potential to improve transportation connectivity and access where it has been infeasible due to, for example, terrain limiting ground transportation infrastructure or the construction of runways for conventional take-off and landing aircraft \cite{gao2023developing}. New AAM services have several potential advantages, including an alternate transportation mode to avoid ground congestion \cite{maheshwari2021evaluating} while providing a clean transportation alternative and environmental benefits \cite{andre2019robust}.

There are significant challenges regarding the implementation, feasibility, and equity implications of introducing a new form of air mobility in complex urban, suburban, and rural airspace. Some of these concerns, such as negative consequences of new infrastructure \cite{takacs2022infrastructural}, directly challenge the aforementioned advantages and benefits of AAM (e.g., climate benefits \cite{zhao2022environmental}, affordability, and accessibility \cite{cohen2023advanced}). Other concerns relate to safely integrating these novel air transport operations into the National Airspace System (NAS), including technical questions regarding air traffic management, demand-capacity balancing, and operational cooperation between AAM fleet operators and civil aviation authorities. The work herein is motivated by a subset of research questions addressing the integration of AAM operations into the NAS. We are motivated by the need to develop rigorous and generalizable models for evaluating new collaborative flow and traffic management roles and strategies within the context of AAM operations.

The concept of strategic air traffic management (ATM) involves ensuring  \emph{demand} (e.g., requested access to constrained resources at airports and within the airspace) and \emph{capacity} (e.g., the nominal or degraded availability of resources) are balanced. While ATM remains an active research area for legacy airspace users, like commercial and general aviation \cite{wandelt2024status}, new ATM paradigms for new airspace users (e.g., AAM aircraft \cite{faa2020uam}, small drones \cite{FAA_UTM_ConOps_v2}, and high-altitude vehicles \cite{lee2024identification}) pose new technical challenges. One example where ATM for legacy and new airspace users may differ, which we tackle in this paper, is the dynamic nature of demand for AAM operations. Whereas commercial airlines have reasonable knowledge of passenger demand weeks (or even months) beforehand, this will not be the case for on-demand air taxis, a potential AAM use case.

\begin{figure}[!htbp]
	\centering
	\includegraphics[width=\linewidth]{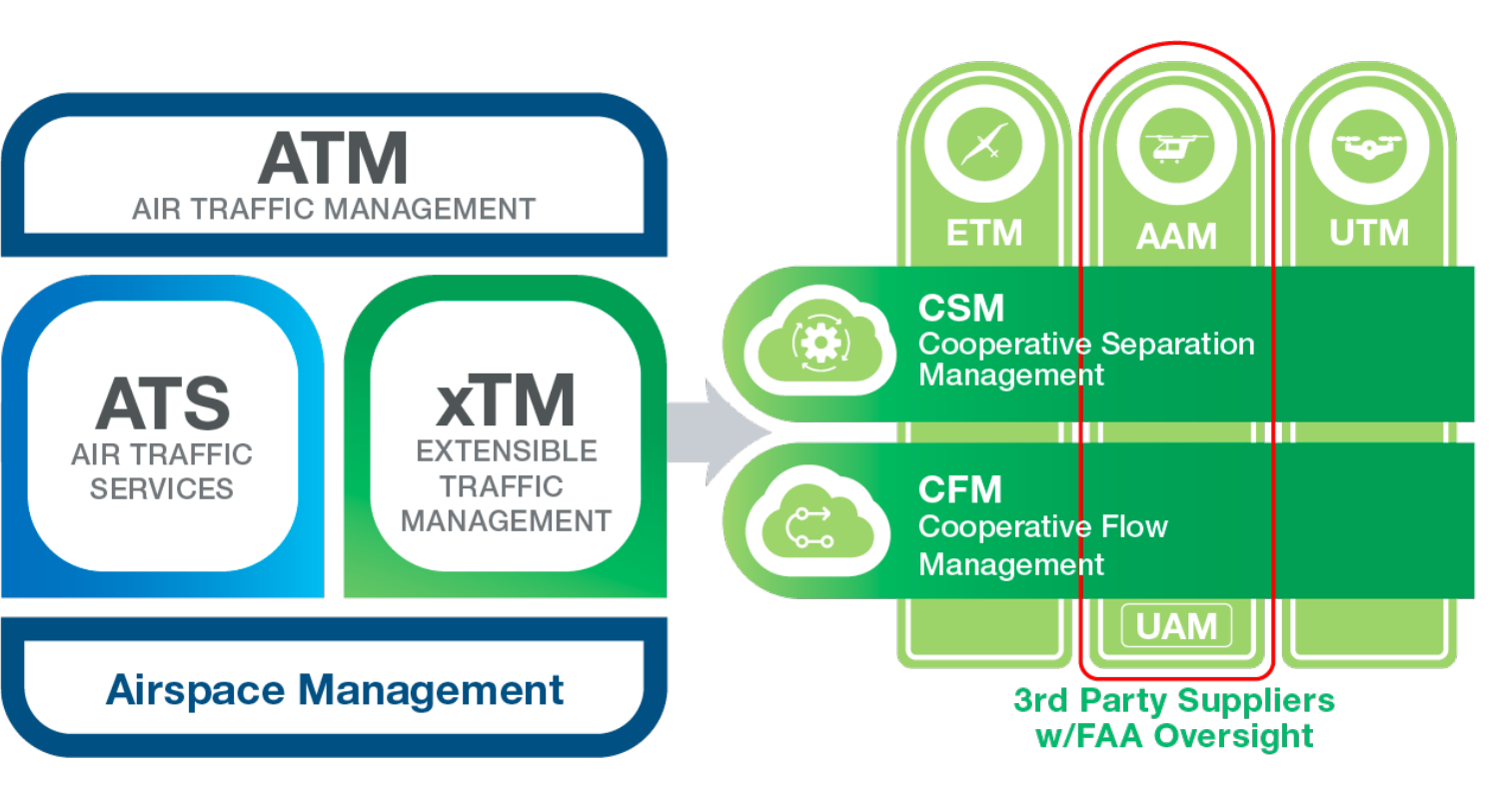}
	\caption{Notional overview of the air traffic, airspace, and traffic management services, modified from the FAA UAM Concept of Operations \cite{FAA_UTM_ConOps_v2}. Air traffic services remain under the purview of the ANSP (FAA in this case), encompassing current air traffic control and traffic flow management within the NAS.}
	\label{fig:faa_utm_conops}
\end{figure}

Given these new challenges, international efforts are underway to define concepts around Uncrewed Aircraft Systems (UAS) Traffic Management (UTM). The initial UTM concept was first developed by NASA in 2013 \cite{kopardekar2014unmanned,kopardekar2016unmanned}, and has now advanced into an implementation plan adopted by the Federal Aviation Administration (FAA) \cite{FAA_UTM_Implementation_Plan_v1_8}. The \emph{Extensible Traffic Management} (xTM) concept, originally defined by NASA and adopted by the FAA, generalizes traffic management concepts beyond UTM and ATM \cite{jung2022overview} to include, e.g., Upper Class E Traffic Management (ETM) for operations above 60,000 feet. Key to the xTM concept is the introduction of \emph{public-private partnership} structures. Traffic management for these new vehicle operations could be performed collaboratively between air navigation service providers (ANSPs), such as the FAA, and fleet operators (e.g., a company providing air taxi services via eVTOLs). Fleet operators may also leverage third-party service providers (e.g., Provider of Services for UAM [PSU] \cite{FAA_UTM_ConOps_v2} and UAS Service Supplier [USS] \cite{rios2019uas}) to manage their operations. Figure \ref{fig:faa_utm_conops} provides an overview of this concept.

Cooperative Operating Practices (COPs), developed by AAM stakeholders and approved by the FAA, will define a set of rules for how AAM fleet operators will conduct their operations \cite{FAA_UTM_ConOps_v2}. The COPs should be designed to ensure safe and efficient operations. In an AAM environment with multiple fleet operators and given potential interactions with other airspace users, the COPs should also ensure fair and equitable operations.

In previous work \cite{wu2025managing}, we modeled the concept depicted in Figure \ref{fig:faa_utm_conops} as \emph{multiple interacting layers}. 
The top layer comprises the ANSP, which provides oversight of the xTM system and protocols, but does not engage directly in day-to-day operations. The middle layer includes one or more service providers responsible for cooperative flow management and critical services enabling safe and efficient operations (e.g., data and information exchange and strategic scheduling). The fleet operator, or multiple fleet operators, reside at the bottom layer. They execute missions which satisfy demand for AAM. We refer readers to \cite{wu2025managing} for additional details on this multi-layered perspective.

\section{Technical Gap and Research Problem}	\label{sec:tech_gap_RP}

Organizing proposed UTM and xTM concepts into these three layers reveals a gap in current modeling literature. There is a plethora of work (e.g., \cite{li2023uav}) related to scheduling and routing of AAM operations by fleet operators (bottom layer). 
There is a reasonable level of characterization for the middle layer containing PSUs and USSs (e.g., \cite{rios2018utm}). The top layer showing ANSP oversight is well-described qualitatively by different aforementioned UTM concepts of operations (e.g., \cite{FAA_UTM_ConOps_v2}). More recently, there have been a number of quantitative modeling approaches that focus on \emph{interactions} between layers, such as between the fleet operators and the airspace service providers (e.g., \cite{chin2022distributed,chin2023protocol}). However, to our best knowledge, until \cite{wu2025managing} there has not been an integrated model encompassing all three layers (ANSPs, airspace service providers, fleet operators). Additionally, \cite{wu2025managing} showed that computational challenges of such a bi-level optimization problem can be overcome with approaches such as neural network-based surrogates.

While \cite{wu2025managing} established a tractable, bi-level optimization approach to the multi-layered xTM concept, it left a number of research directions worth exploring: incorporating more responsive and dynamic AAM demand models, taking into consideration different congestion management strategies, and factoring in pop-up AAM flights. Tackling these limitations from \cite{wu2025managing} is further motivated by the identified technical barriers towards broader AAM adoption, namely fleet and airspace management as well as lack of demand-centric considerations \cite{patterson2021initial,garrow2025market}. Complementing the other UTM models which we will survey in Section \ref{sec:lit_review}, we envision the work in this paper contributing to a set of rigorous modeling and simulation tools that may inform the integration of AAM operations into the NAS (or airspace managed by other ANSPs\footnote{We note that although we focus on the US-centric UTM vision, analogous concepts such as U-Space in Europe have been developed \cite{Riccardi2024}. 
}). Furthermore, the analysis described herein can inform the development of COPs for demand-capacity balancing and congestion management approaches as well as the use of simulation models for evaluating COPs in AAM environments. 



Under the xTM concept where day-to-day operations are managed between the fleet operators and the airspace service providers, with infrequent oversight actions taken by the ANSP, we seek to answer the following questions:

\begin{enumerate}
	\item How does time-varying demand impact airspace congestion and delays, and what kind of mitigating actions (e.g., congestion management) can be taken?
	\item How does the inclusion of AAM flights that are unscheduled (perhaps due to being of a much higher priority, such as an air ambulance flight) impact aforementioned mitigating actions?
\end{enumerate}

\noindent
To answer these questions, we adopt and extend the bi-level optimization formulation from \cite{wu2025managing}. We now give an  overview of the technical approach, which we mathematically detail in Section \ref{sec:model_formulation}. We model the interactions between the bottom layer (fleet operators) and middle layer (airspace service providers) as a scheduling optimization problem---the \emph{low-level planner} (LLP)---over a network of Origin-Destination (OD) pairs. This interaction between the bottom and middle layers are monitored by the ANSP. We model this as a second optimization problem---the \emph{high-level planner} (HLP)---that adjusts network path costs and vertiport\footnote{A \emph{vertiport} to AAM is analogous to an airport for legacy air transportation operations. There are several different facets related to vertiports which are also active areas of research (see, e.g., \cite{schweiger2022urban}). For simplicity, we assume all AAM trips in our model begin and end at some generic vertiport within our AAM network.} landing costs to achieve a given congestion management objective. This is a form of \emph{congestion pricing}. 

\begin{figure*}[!htbp]
	\centering
	\includegraphics[width=0.75\linewidth]{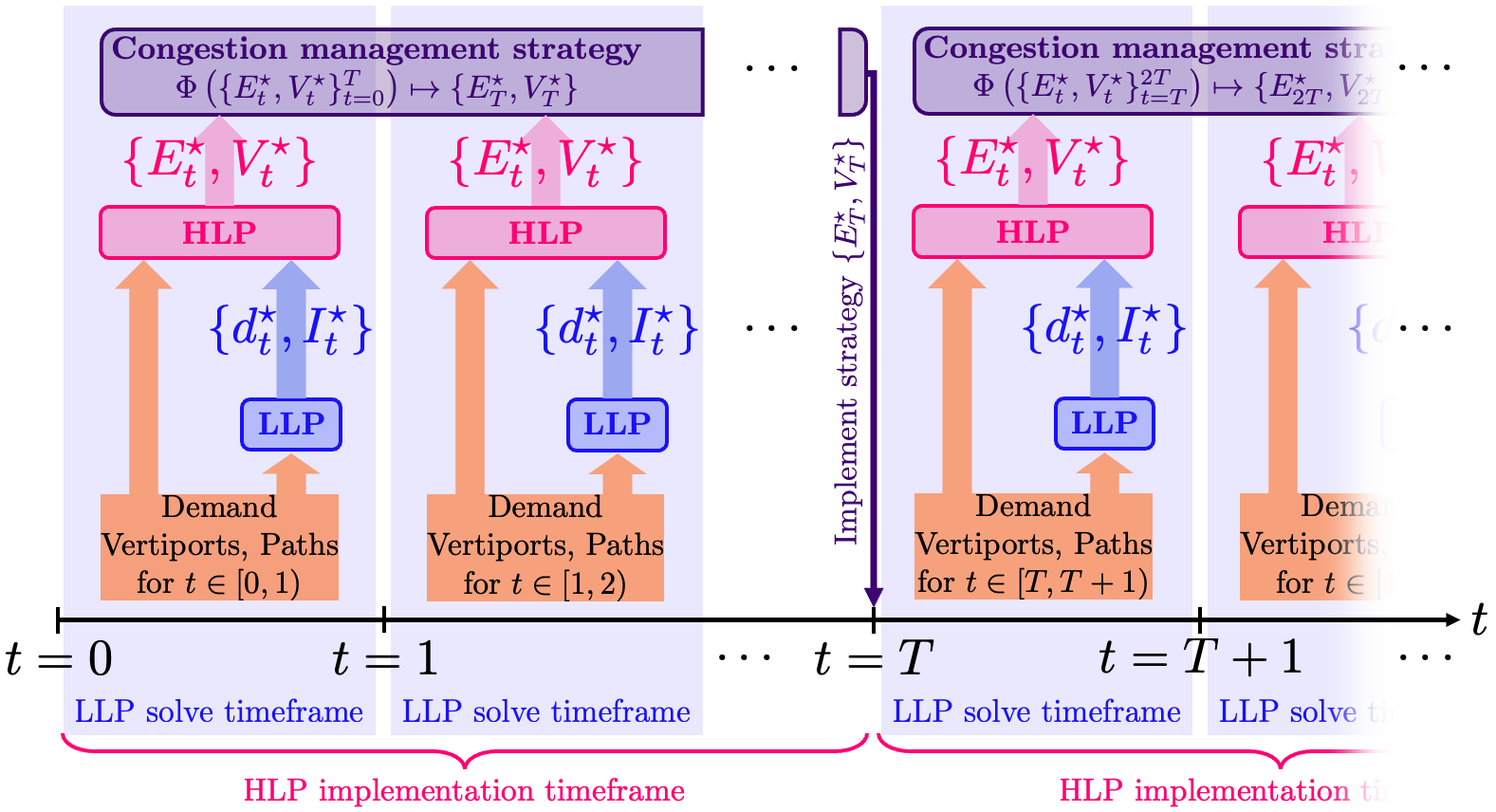}
	\caption{Overview of the model workflow for the AAM congestion management bi-level optimization problem.}
	\label{fig:ov1}
\end{figure*}

Figure \ref{fig:ov1} provides a graphical depiction of the interleaved LLP and HLP. As the LLP mimics day-to-day operations, its solution is implemented on a much shorter time horizon compared to the HLP. The HLP, which mimics an oversight role, considers congestion outcomes at each LLP solve timeframe; the overall congestion management strategy suggested by the HLP is implemented after multiple LLP solve timeframes. The two optimization problems are coupled as a bi-level optimization problem, drawn from the structure of Stackelberg games, addressed in Section \ref{sec:model_formulation}. 
We note that the time indices in Figure \ref{fig:ov1} are notional; we detail the LLP and HLP timescales used in the experimental design in Section \ref{sec:exp_setup_results}. The mathematical notations used in Figure \ref{fig:ov1} are also introduced in the LLP and HLP formulations in Section \ref{sec:model_formulation}.


%



\section{Literature review}	\label{sec:lit_review}

We focus this literature review on congestion pricing and bi-level optimization models within the aviation context, UTM and xTM, as well as demand estimation for AAM. For a more comprehensive literature review on AAM topics, we refer readers to, for example, \cite{garrow2021urban}.

\subsection{Congestion pricing and bi-level optimization models}

Congestion pricing as an approach for managing a network of small drones was explored in \cite{wang2023learning}, although a pure learning-based algorithm was used to derive pricing policies. Such learning-based approaches lead to difficulties in understanding why specific pricing actions were taken. Within air transportation, the use of congestion pricing has typically focused on airports, relating to airport slot pricing  \cite{pels2004economics,verhoef2010congestion}. This perspective was extended to air traffic flow management and \enquote{trading} delays in \cite{castelli2011design}, as well as congested air traffic networks \cite{adler2022competition}. However, these are completely different contexts compared to the setting for UTM and xTM. Finally, \cite{he2024distributed} examines route planning while considering congestion pricing, but does not focus on strategic traffic flow coordination or management---it also does not consider the public-private partnership structure envisioned for UTM and xTM. 

Given the nested structure of the LLP and HLP (i.e., multiple runs of the LLP exist within an implementation period of the HLP), we formulated this problem as a \emph{bi-level optimization}; \cite{sinha2017review} provides an in-depth review of bi-level optimization problems. Furthermore, we cast this as a Stackelberg game \cite{frantsev2012finding}, where the ANSP (through the HLP) is the \emph{leader} attempting to achieve a specific congestion management goal across the entire AAM network. The airspace service providers are \emph{followers}, and optimize schedules and route options (via the LLP), interacting with the fleet operators. Such Stackelberg games can be solved via bi-level optimization, e.g., see \cite{yue2017stackelberg}.

\subsection{UTM and xTM}

At the level of scheduling and assigning delays to individual UAS and AAM vehicles, \cite{chin2023protocol} examines protocol-based approaches such as via back-pressure metrics. In \cite{chin2021efficiency}, strategic traffic flow management for UAS was modeled via integer programming, analyzing the trade-off between efficiency and fairness. Additional studies on fairness in UTM include \cite{evans2020fairness}, where 
the authors showed delay imbalances can occur between fleet operators dependent on when flight plans were filed. Technical challenges imposed by centralized (e.g., one airspace with one PSU) versus decentralized (e.g., multiple airspace regions, multiple PSUs) settings were studied in \cite{chin2022distributed} and \cite{qin2023market}, with the former focusing on distributed optimization and the latter focusing on mechanism design (e.g., auctions). Finally, \cite{skorup2019auctioning} examines airspace auction concepts from a regulatory perspective, but does not connect this back to specific strategic traffic management actions.

\subsection{AAM demand estimation and modeling}

Although we do not consider any specific use cases in this paper, AAM demand estimation has been examined for applications such as airport shuttles \cite{roy2021user,roy2022flight}. Additional realism may be added to our demand function models by incorporating passenger numbers estimated through studies such as \cite{haan2021commuter} or \cite{boddupalli2024mode}. Our use of different demand profiles follows from similar concepts in commercial air transport, such as the Boeing Decision Window Model, which has been used to reason about aggregate passenger demand behavior \cite{parker2007passengers,berge2007future}.

\section{Contributions of Work}	\label{sec:contributions}

This paper builds significantly on the bi-level optimization model from \cite{wu2025managing} by introducing and studying the impacts of several extensions to the original model. These extensions include parameterized demand functions, comparing different congestion management strategies, and stochastic, \emph{unscheduled} demand. Specifically, our contributions in this paper are as follows:

\begin{enumerate}
	\item We generate dynamic, time-varying demand functions to better simulate OD trips within our AAM network. We parameterize them to allow for ease of building different demand scenarios (e.g., morning versus afternoon peaks). 
	\item We explore different congestion management strategies for adjusting the per-path tolling strategies and vertiport landing fees. We compare the performance of these strategies for different demand patterns.
	\item We study the impact of unscheduled \enquote{pop-up} AAM flights with varying likelihood of occurrence. We study the degradation of the congestion management strategies given the pop-up flights.
\end{enumerate}

%
%
%
%
%

\section{Model Formulation}	\label{sec:model_formulation}
We begin by defining notations that will be used throughout the remainder of the paper: We denote by $T$ the planning horizon (i.e., the LLP solution timeframe in Figure \ref{fig:ov1}), by $A$ the set of airspace sectors, by $P$ the set of paths (e.g., AAM corridors), by $S$ the set of PSUs, by $V$ the set of vertiports considered in our network, and by $F$ the set of all flights. Refinements on the set of flights include $F(s, \tau)$ the set of flights governed by PSU $s$ at time $\tau$, $F^{(pp)}(s,\tau)$ the set of unscheduled (i.e., pop-up) flights governed by PSU $s$ at time $\tau$, $F(a)$ the set of flights that will possibly fly through airspace $a$, $F(v)$ the set of flights $f$ that will land at vertiport $v$, and $FP(a)$ the set of flight and path pairs where flight $f$ will fly through airspace $a$ along path $p$.

We now define notation pertaining to the LLP and HLP functionalities. We denote by $C^{d}_{f}$ the unit delay cost for flight $f$, by $C^{p}_{p}$ the path cost for path $p$, by $C^{v}_{v,\tau}$ the vertiport landing fee for vertiport $v$ at time $\tau$. The specific scheduled departure time for flight $f$ is denoted by $t^{d}_{f}$, the estimated flight time along path $p$ (i.e., en route or airborne time) is denoted by $t^{p}_{p}$, and the actual departure time of flight $f$ is denoted by $d_{f}$. Since we formulate the LLP as an integer program, we define $I_{f,p}$ as a binary variable indicating whether flight $f$ will choose path $p$, and $z_{f,\tau}$ is a binary variable indicating whether flight $f$ will arrive at the destination at time $\tau$. We define  $\lambda_{f,p,\tau}$ as the binary variable that indicates whether flight $(f,p)\in FP(a), \forall a \in A$ arrives at airspace $a$ during time period $\tau$ along path $p$, by $\kappa_{f,\tau}$ the binary variable that specifies whether flight $f \in F(a), \forall a \in A$ will pass through airspace $a$ during time period $\tau$, by $\mu_{a,\tau}$ the binary variable that indicates whether the traffic count at airspace $a$ during time period $\tau$ exceeds its capacity $Q^{A}_{a}$, and by $\Delta_{a,\tau}$ the auxiliary variable that indicates the congestion at airspace $a$ during time period $\tau$. We will use $\epsilon^{p}_{p}$ to denote the \emph{tolling strategy} on path $p$, and $\epsilon^{v}_{v,\tau}$ the \emph{landing fee strategy} at vertiport $v$ at time $\tau$. Finally, we denote by $Q^{Dep}_{v}$ the departure capacity of vertiport $v$, and by $Q^{Arr}_{v}$ the arrival capacity of vertiport $v$.

\subsection{Low-level planner (LLP)} \label{subsec:LLP}
As mentioned above, the LLP determines optimal departure times and path selections for a fleet operator so that costs will be minimized. The formulation of the LLP is as follow:
\begin{subequations}
	\begin{alignat}{2}
		\min_{d, I} & \sum_{\tau \in T} \sum_{s \in S} \sum_{f \in F(s,\tau)} \left( C^{d}_{f}(d_{f} - t^{d}_{f}) + \sum_{p \in P} I_{f,p} C^{p}_{p} \right)  \nonumber \\ 
		& + \sum_{v \in V}\sum_{f\in F(v)}\sum_{\tau \in T} C^{v}_{v,\tau} z_{f,\tau} \label{eq:low_level_obj} \\
		\text{s.t.} \quad & \|d_{f} - d_{f'}\| \geq \frac{3600}{Q^{\text{Dep}}_{v}}, \quad \forall \tau \in T, v \in V, s \in S, s' \in S, \nonumber \\ 
		& f \in F(s,\tau), f' \in F(s',\tau), \label{eq:low_level_const_1}\\
		& \|d_{f} + \sum_{p \in P} I_{f,p} t^{p}_{p} - d_{f'} - \sum_{p \in P} I_{f',p} t^{p}_{p}\| \geq \frac{3600}{Q^{\text{Arr}}_{v}}, \nonumber \\ 
		& \forall \tau \in T, v \in V, s, s'\in S, f \in F(s,\tau), f' \in F(s',\tau), \label{eq:low_level_const_2}\\
		& d_{f} + \sum_{p \in P} I_{f,p} t^{p}_{p} \geq \tau - M(1 - z_{f,\tau}), \quad \forall f \in F, \tau \in T, \label{eq:low_level_const_3}\\
		& d_{f} + \sum_{p \in P} I_{f,p} t^{p}_{p} \leq \tau + 1 + M(1 - z_{f,\tau}), \forall f \in F, \tau \in T, \label{eq:low_level_const_4}\\
		& d_{f} - t^{d}_{f} \geq 0, \quad \forall f \in F, \label{eq:low_level_const_5}\\
		& \sum_{p \in P} I_{f,p} = 1, \quad \forall f \in F, \label{eq:low_level_const_6}\\
		& \sum_{\tau \in T} z_{f,\tau} = 1, \quad \forall f \in F. \label{eq:low_level_const_7}
	\end{alignat}
	\label{eq:low_level}
\end{subequations}
The objective of the LLP \eqref{eq:low_level_obj} is to minimize the total cost, which includes the delay cost, path selection cost, and vertiport landing fees. Constraints \eqref{eq:low_level_const_1}–\eqref{eq:low_level_const_2} define the vertiport departure and arrival requirements, ensuring sufficient time separations between consecutive landings and departures. Constraints \eqref{eq:low_level_const_3}–\eqref{eq:low_level_const_4} account for the time-dependent nature of vertiport landing fees, modeling the varying costs at different times of the day. Specifically, these constraints determine the time when a flight lands. Lastly, constraints \eqref{eq:low_level_const_5}–\eqref{eq:low_level_const_7} ensure that no flights depart early and that each flight has selected exactly one path and one landing time.

\subsection{High-level planner (HLP)} \label{subsec:HLP}

Recall from Section \ref{sec:tech_gap_RP} that we motivated the formulation of the relationship between the HLP and the LLP as a bi-level optimization problem. The HLP has an overarching objective of congestion management across the entire network, whereas the LLP focuses on how fleet operators might schedule and route AAM flights, subject to delay costs and tolls. The HLP decisions impact the LLP by imposing different path- and vertiport-based tolls. The LLP decisions impact the HLP by changing the level of congestion based on the interplay between fleet operator demand and capacities. The HLP (with the LLP nested within it) is formulated as follows:
\begin{subequations}
	\begin{alignat}{2}
		&\min_{\epsilon^{p}, \epsilon^{v}, I, d, \Delta, \kappa, \lambda, \mu} \qquad
		\sum_{a \in A} \sum_{\tau \in T} \Delta_{a, \tau} \label{eq:high_level_obj} \\
		&\text{s.t.} \quad \nonumber \\
		& (I^{\star}, d^{\star}) = \underset{I, D}{\argmin} \sum_{\tau \in T} \sum_{s \in S} \sum_{f \in F(s, \tau)} \Bigl(C^{d}_{f} - t^{d}_{f}  \nonumber \\ 
		&+ \sum_{p \in P} I_{f, p}(C^{p}_{p} + \epsilon^{p}_{p})\Bigr) + \sum_{v \in V} \sum_{f \in F(v)} \sum_{\tau \in T}(C^{v}_{v, \tau} + \epsilon^{v}_{v, \tau})z_{f, \tau}, \label{eq:high_level_c1} \\
		& d_{f} + t^{p}_{p} \geq \tau - M(1 - \lambda_{f, p,\tau}), \; \forall (f,p) \in FP(a), \forall \tau \in T, \label{eq:high_level_c2}\\
		& d_{f} + t^{p}_{p} \leq \tau + 1 - M(1 - \lambda_{f, p,\tau}), \forall (f,p) \in FP(a), \forall \tau \in T, \label{eq:high_level_c3} \\
		& \sum_{\tau \in T} \lambda_{f,p,\tau} = 1, \quad \forall (f,p) \in FP(a), \label{eq:high_level_c4} \\
		& \kappa_{f, \tau} \leq I_{f, p}, \quad \forall \tau \in T, \forall (f, p) \in FP(a), \forall a \in A, \label{eq:high_level_c5} \\
		& \kappa_{f, \tau} \leq \lambda_{f,p,\tau}, \quad \forall (f,p) \in FP(a), \forall \tau \in T, \label{eq:high_level_c6} \\
		& I_{f, p} + \lambda_{f, \tau} - 1 \leq \kappa_{f, \tau}, \quad \forall \tau \in T, \forall (f, p) \in FP(a), \forall a \in A, \label{eq:high_level_c7} \\
& \Delta_{a, \tau} \leq \bigl(\sum_{f \in F(a)} \kappa_{f, \tau} - Q^{A}_{a}\bigr) + M(1 - \mu_{a, \tau}),  \label{eq:high_level_c8} \\
		& \bigl(\sum_{f \in F(a)} \kappa_{f, \tau} - Q^{A}_{a}\bigr) \leq \Delta_{a, \tau}, \quad \forall a \in A, \forall \tau \in T, \label{eq:high_level_c9} \\
		& 0 \leq \Delta_{a, \tau}, \quad \forall a \in A, \forall \tau \in T, \label{eq:high_level_c10} \\
		& \Delta_{a, \tau} \leq M \mu_{a, \tau}, \quad \forall a \in A, \forall \tau \in T, \label{eq:high_level_c11} \\
		& \text{Low-level planner (LLP) constraints \eqref{eq:low_level_const_1}-\eqref{eq:low_level_const_7}.} \nonumber
	\end{alignat}
	\label{eq:high_level_1}
\end{subequations} 
The optimization problem in \eqref{eq:high_level_1} couples the decision-making between the HLP and the LLP. The optimization problem's structure follows the general form a bi-level optimization problem, given explicitly in Definition \ref{bi-level program}:

\begin{defn}[Bi-level program \cite{sinha2017review}]
	For the upper-level objective function $F$ and lower-level objective function $f$, a \emph{bi-level program} is given by
	\begin{equation}
			\begin{aligned}
					&\min_{x_{u}\in X_{U}, x_{l}\in X_{L}}\quad F(x_u,x_l)\\
					\mathrm{s. t. } \; & x^{\star}_l \in \underset{x_l \in X_L}{\argmin} \{f(x_{u},x_{l}):g_{j}(x_u,x_l) \leq 0, j = 1, \dots, J\}, \\
					&G_{k}(x_u,x_l) \leq 0, k = 1, \dots, K, 
				\end{aligned}
		\end{equation}
	where $X_{U}, X_{L}$ are sets of decision variables for the upper- and lower-level problems, respectively. $G_{k}$ with $k = 1, \dots, K$ denote the upper-level constraints, and $g_{j}$ with $j = 1, \dots, J$ represent the lower-level constraints, respectively. 
	\label{bi-level program}
\end{defn}

	The objective of the HLP is to reduce congestion in each airspace volume $a$ during each time period $\tau$, represented by $\Delta_{a,\tau} \coloneqq \max(0, \sum_{f \in F} I_{f, \tau} - Q^{A}_a)$. Note that we have defined a \emph{congestion metric} $\Delta_{a, \tau}$ which counts the number of AAM flights (i.e., the traffic count) which exceed the airspace capacity.
Constraint \eqref{eq:high_level_c1} corresponds to the LLP optimality condition in \eqref{eq:low_level}, ensuring that $I^{\star}$ and $d^{\star}$ are optimal solutions of the LLP under the imposed high-level strategies (i.e., path and vertiport tolls). Constraints \eqref{eq:high_level_c2}–\eqref{eq:high_level_c4} determine the time when a flight $f$ traverses airspace $a$, while constraints \eqref{eq:high_level_c5}–\eqref{eq:high_level_c7} calculate the number of flights passing through airspace $a$ during time period $\tau$. 
To handle the maximization function within the objective $\Delta_{a, \tau}$, we introduce an auxiliary variable $\Delta_{a, \tau}$, representing the congestion in airspace $a$ during time period $\tau$. Constraints \eqref{eq:high_level_c8}–\eqref{eq:high_level_c11} implement the maximization function using the big-$M$ method.

The proposed bi-level program in \eqref{eq:high_level_1} provides a complete and rigorous description of the HLP and its interactions with the LLP. However, there are significant tractability challenges when it comes to solving this problem via the current formulation: Constraint \eqref{eq:high_level_c1} necessitates that $(I^{\star}, d^{\star})$ be the optimal solution to the LLP. The objective function of the LLP depends directly on the imposed tolling strategies from the HLP, and there is no guarantee that the objective function of the LLP is convex with respect to the HLP congestion management strategies. Hence, the bi-level program in \eqref{eq:high_level_1} is intractable, and not amenable to single-level reduction techniques (e.g., \cite{sinha2017review}) for solving analogous bi-level programs. To address this challenge, in previous work we proposed a Neural Network (NN)-based surrogate approach to approximate the LLP objective function (value function approximation) \cite{wu2025managing}. For further details on the construction of the NN surrogate and the derivation of the tractable formulation of the HLP, we refer readers to \cite{wu2025managing}. 

\subsection{Congestion management strategies via aggregation} \label{subsec:aggregation}
As shown in Figure \eqref{fig:ov1}, after each execution of the HLP, we obtain a set of path and vertiport tolling decisions. By the end of each HLP implementation timeframe, we need to \emph{aggregate} all HLP decisions generated within that timeframe and propose a single congestion management strategy for operations in subsequent time periods. We propose two aggregation strategies, henceforth known as congestion management strategies. The first strategy is a naive averaging approach over the same time periods of the day. This method assumes AAM demands during the same time periods on different days exhibit similar patterns. 

The second type of strategy is based on the assumption that demand patterns may vary across different days. Therefore, we aim to develop a strategy that captures the similarity between historical demands and future demands. To achieve this, we utilize the Wasserstein distance \cite{villani2009wasserstein} to measure the distance between demand patterns. We then derive inverse distance-based weights to aggregate historical HLP decisions. The inverse distance is important as more similarity translates to smaller Wasserstein distances, and vice versa. The second strategy is implemented as follows: let $\mathrm{d}^{n}_{m}$ represent the AAM demand for day $n$ and planning horizon $m$. To apply the Wasserstein distance between two distributions, it is necessary to ensure they share the same support and that the total measure on the support is equal to 1. Consequently, we normalize the demand and represent the distribution as
\begin{equation}
	\begin{aligned}
		\mathbf{T}_{m} = \frac{\mathrm{d}_{m}}{\sum_{m}\mathrm{d}_{m}}, \forall m = 1,2,\dots, M.
	\end{aligned}
\end{equation}
Therefore, the Wasserstein distance between the demand patterns of the future day and one of the historical days is given by $\mathbf{d}_{W}(\widehat{\mathbf{T}}, \mathbf{T}^{(n)})$, where we define $\widehat{\mathbf{T}}$ as the distribution for the test day and $\mathbf{T}^{(n)}$ as the distribution for a historical day. Based on the Wasserstein distance between demand patterns, we define the inverse weights for each historical day as 
\begin{equation}
	\begin{aligned}
		\mathbf{w}_{n} = \frac{ \left( \mathbf{d}_{W}(\widehat{\mathbf{T}}, \mathbf{T}^{(n)}) + \mathbb{\varepsilon} \right)^{-1} }{\sum_{n=1}^{N} \left(  \mathbf{d}_{W}(\widehat{\mathbf{T}}, \mathbf{T}^{(n)}) + \mathbb{\varepsilon}\right)^{-1} },
	\end{aligned}
\end{equation}
where $\varepsilon > 0$ is a very small, positive scalar to avoid dividing by zero. With the proposed weights, denote by $\epsilon^p_m, \epsilon^v_m$ the second congestion management strategy to be imposed; these are calculated as follows:
\begin{equation}
	\begin{aligned}
		\epsilon^{p}_{m} =\sum_{n} \mathbf{w}_{n} \epsilon^{p^{(n)}}_{m}; \; \epsilon^{v}_{m} = \sum_{n}\mathbf{w}_{n} \epsilon^{v^{(n)}}_{m}, \quad \forall m = 1, \dots, M.
	\end{aligned}
\end{equation}
\subsection{Modeling of unscheduled, pop-up AAM flights}
To model the integration of unscheduled flights with the current scheduled flights, we assume all unscheduled flights are given higher priority by the LLP (e.g., air ambulances, emergency response drones). This prioritization ensures that the LLP always assigns the shortest path to unscheduled flights, while other flights must be rescheduled to satisfy vertiport capacity constraints. In this scenario, on-demand flight schedules are treated as hard constraints, and these constraints are incorporated into the LLP:
\begin{equation}
	\begin{aligned}
		& \|d_{f} - d_{f'}\| \geq \frac{3600}{Q^{\text{Dep}}_{v}}, \\
		&\quad \forall \tau \in T, \, v \in V, \, (s,s') \in S \, f \in F(s,\tau), \, f' \in F^{(pp)}(s, \tau),  \\
		& \|d_{f} + \sum_{p \in P} I_{f,p} t^{p}_{p} - d_{f'} - \sum_{p \in P} I_{f',p} t^{p}_{p}\| \geq \frac{3600}{Q^{\text{Arr}}_{v}}, \\
		&\quad \forall \tau \in T, \, v \in V, \, (s,s') \in S, \, f \in F(s,\tau), \, f' \in F^{(pp)}(s, \tau).
	\end{aligned}
\end{equation}

\section{Experimental Setup and Results}	\label{sec:exp_setup_results}

\subsection{Demand function} \label{subsec:demand}
\begin{figure}[th]
	\centering
	\begin{subfigure}[b]{0.4\textwidth} 
		\centering
		\includegraphics[width=\textwidth]{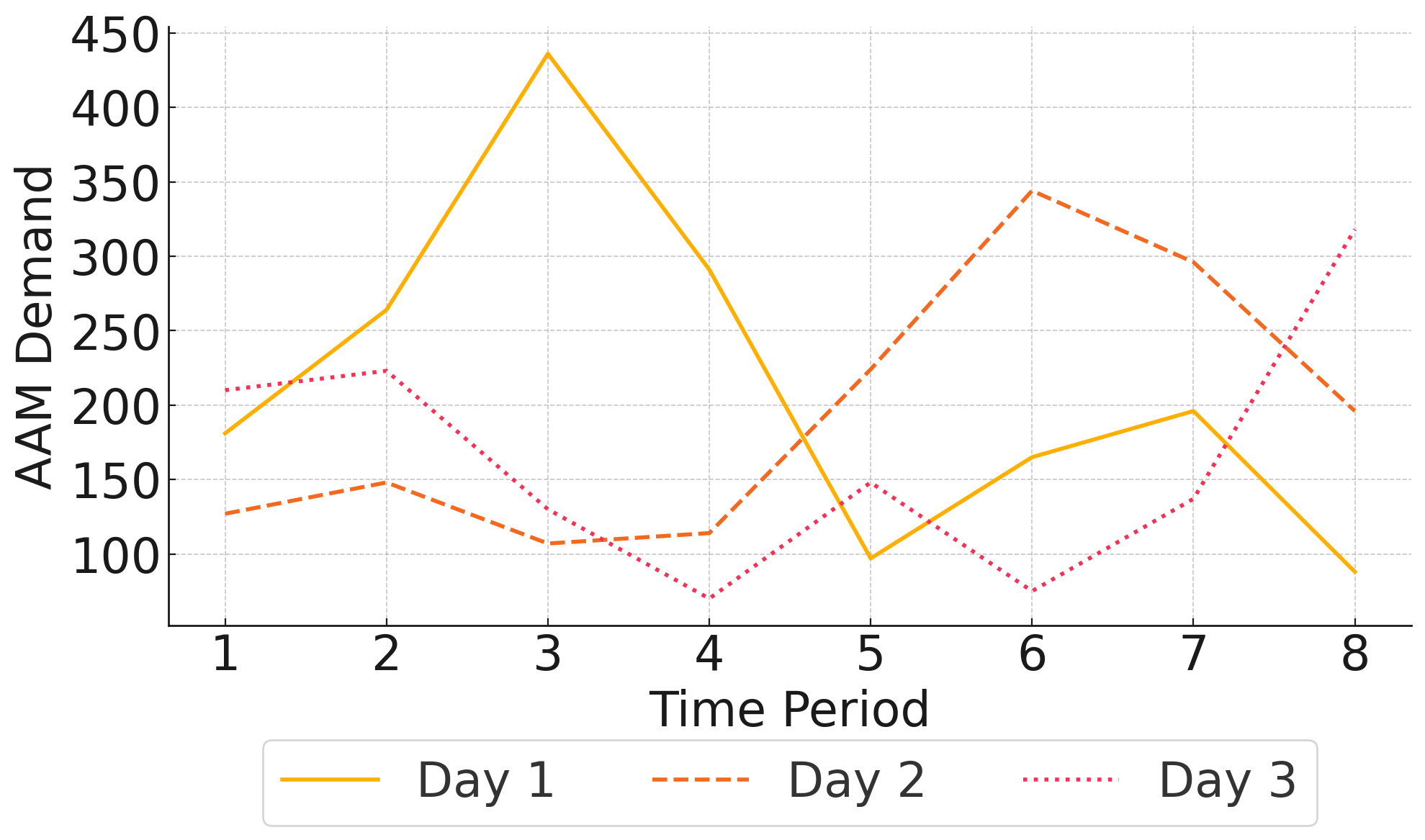}
		\caption{Demand profiles for test set 1.}
		\label{fig:AAM_demand_test1}
	\end{subfigure}
	\hfill
	\begin{subfigure}[b]{0.4\textwidth} 
		\centering
		\includegraphics[width=\textwidth]{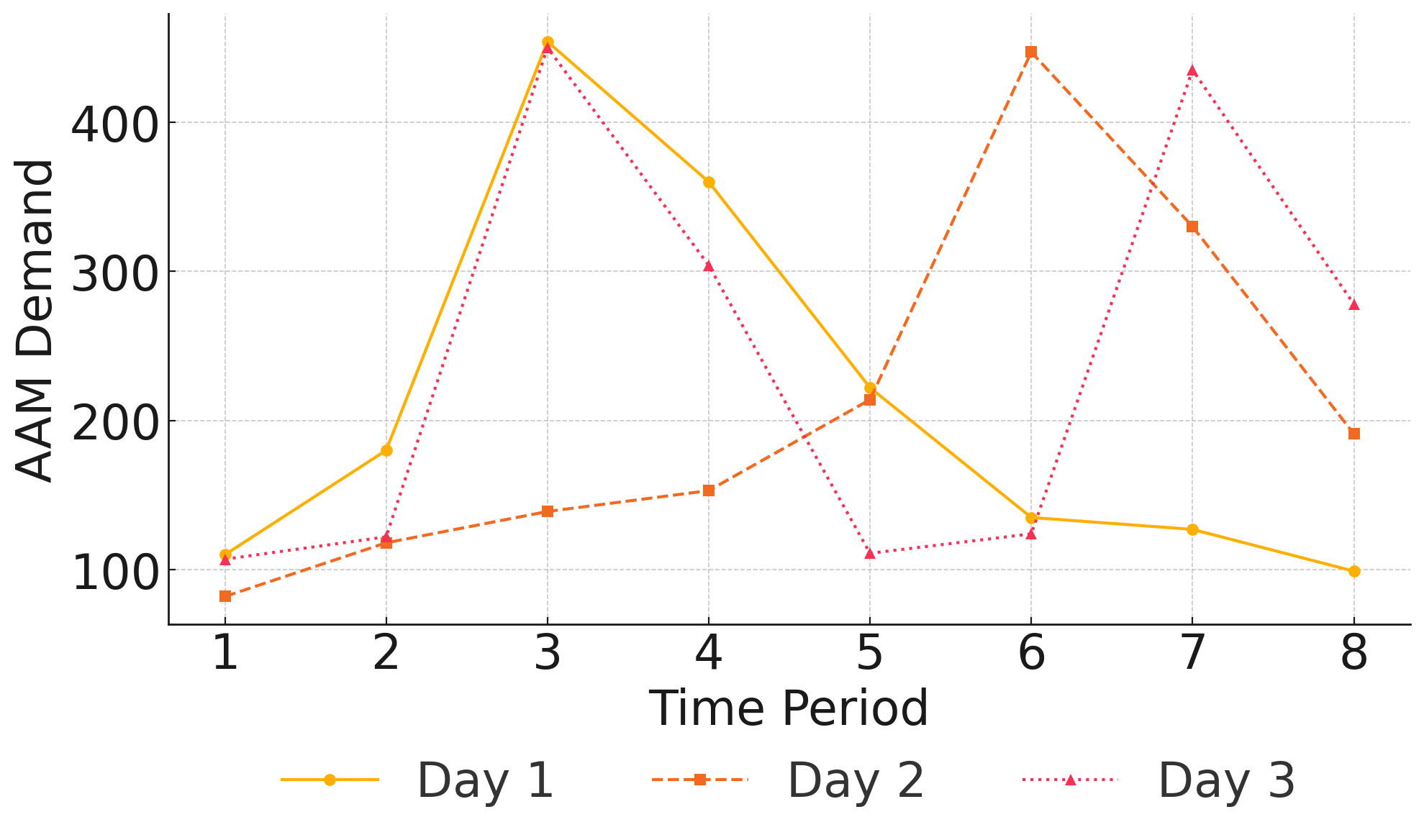}
		\caption{Demand profiles for test set 2.}
		\label{fig:AAM_demand_test2}
	\end{subfigure}
	\caption{AAM demand profiles the two test sets, across 3 days with each day split into 8 different 3-hour time periods.}
	\label{fig:AAM_demand_comparison}
\end{figure}

We first describe the underlying demand functions used to generate flight schedules. We chose a stylized demand function as our focus is on congestion management, but this demand function could be swapped out by data-driven ones such as those derived in \cite{roy2021user}. To generalize the demand function, we assume AAM demands for the OD pair $(i, j) \in V \times V$ are parameterized via a sinusoidal function:
\begin{equation}
	\begin{aligned}
		\mathbf{D}_{i,j}(t) = \mathbf{A}_{i,j}(t)\sin(\mathbf{B}t + \mathbf{C}_{i,j}) + \Gamma_{i,j}(t) + \mathbf{e}_{i,j}(t),
	\end{aligned}
\end{equation}
where $\mathbf{A}$, $\mathbf{B}$, $\mathbf{C}$, and $\Gamma$ represent the time-varying amplitude, frequency, phase shift, and baseline of the sinusoidal function, respectively. We also add a Gaussian noise term $\mathbf{e}_{i,j} \sim \mathcal{N}(0,\sigma)$ to model demand stochasticity.\footnote{Stochasticity here refers to random, per-OD increases (or decreases) in the number of AAM flights enabled by the $\mathbf{e}_{i,j}$ term. We slightly overload the term \emph{demand stochasticity} by also using it when referring to the presence of unscheduled, pop-up flights (see Section \ref{ssec:unscheduled_demand} for more details).} By adjusting the frequency and phase shift, we control the number and timing of peak hours, while the amplitude and baseline determine the magnitude of demands.  The time-varying demand $\mathbf{D}_{i,j}(t)$ is treated as the mean value of a Poisson distribution, from which interarrival times can be sampled via an exponential distribution, i.e., $\mathcal{T}_{i,j}(t) \sim \mathrm{Exp}\left(1 / \mathbf{D}_{i,j}(t)  \right)$, thereby generating the desired flight schedules.
Additionally, we classify all vertiports as either \emph{busy} or \emph{non-busy}, with different levels of demand. We consider three types of demand functions corresponding to scenarios of morning, afternoon, and evening peaks. To further demonstrate the generalizability of the HLP congestion management strategies under stochastic demands, we generate two distinct sets of demand profiles with differing amplitudes, baselines, and timings of peaks. These demand sets are shown in Figures \ref{fig:AAM_demand_test1} and \ref{fig:AAM_demand_test2}.

\subsection{Neural Network surrogate model: Setup and performance}

We now discuss the dimensionality of inputs and the stability of the NN surrogates, an aspect not examined in \cite{wu2025managing}. In our numerical experiments, we consider a vertiport network based in Los Angeles, California, consisting of seven vertiports and a total of 126 paths. For each planning horizon, if we allow path tolls on all paths and hourly vertiport landing fee strategies at all vertiports, the input dimension of the NN becomes $|T| \times 7 + 126$. In order to train the NN surrogate, we must generate training samples from solving the LLP. However, due to the computational complexity of the LLP, we found it logistically feasible to generate only 1,000 LLP samples as training data for the LLP solve timeframe. 

To evaluate the performance of the NN surrogate, for an arbitrary LLP solve timeframe, we denote by $\epsilon^{p^{\star}}$ and $\epsilon^{v^{\star}}$ the tolling decisions generated by the HLP, $\phi(\epsilon^{p^{\star}}, \epsilon^{v^{\star}})$ as the true optimal value of the LLP, and $\widehat{\phi}(\epsilon^{p^{\star}}, \epsilon^{v^{\star}})$ as the NN surrogate approximation of the optimal value. We can now define a performance metric for the NN surrogate, the \emph{approximation ratio}, which is given by: 
\begin{equation}
	\begin{aligned}
		\textbf{Approx. Ratio} = \frac{|\widehat{\phi}(\epsilon^{p^{\star}}, \epsilon^{v^{\star}}) - \phi(\epsilon^{p^{\star}}, \epsilon^{v^{\star}})|}{\phi(\epsilon^{p^{\star}}, \epsilon^{v^{\star}})}.
	\end{aligned}
\end{equation}
If the approximation ratio is 0, then the NN surrogate perfectly matches the analytical optimal solution of the LLP. Thus, lower approximation ratios are desired. To analyze the impact of the input dimension on the NN surrogate's performance, we select an arbitrary day with eight 3-hour LLP solve timeframes and compare two different input sets. The input dimension of the first set is 147, while the input dimension of the second set is 8. As shown in Tables \ref{tab:Approx_ratio_full} and \ref{tab:Approx_ratio_reduced}, the first set of NN surrogates performs poorly, with approximation ratios as large as almost $43\%$. From our experimental results, we observed that high approximation ratios led to unstable and non-reproducible high-level decisions. This was a shortcoming in \cite{wu2025managing}, where the input dimensions of the NN surrogates were very large.

In contrast, we can limit the set of paths and vertiports where tolling and landing fees are imposed. This is operationally sensible. For example, in surface congestion pricing and tolling, only major, congestion-prone arterial roads are tolled \cite{NAP23427}, and not every single street within the road network. If we limit the set of paths and vertiports to only highly-congested ones, the approximation ratios improve significantly, ranging from $0.11\%$ to approximately $5\%$. Consequently, in this paper, rather than controlling every path and vertiport in the system, we allow for only a subset of paths and vertiports to be tolled to achieve better NN surrogate performance.
\begin{table}[h!]
	\centering
	\caption{Optimal objective function values and approximation ratios with NN surrogate input dimension of 147.}
	\begin{tabular}{|c|c|c|}
		\hline
		\textbf{Obj. Val (\$)} & \textbf{Approx. Val (\$)} & \textbf{Approx. Ratio ( \%)} \\
		\hline
		5252 & 3208 & 38.92 \\
		6403 & 3724 & 41.84 \\
		4444 & 2777 & 37.50 \\
		5458 & 3124 & 42.76 \\
		10190 & 6167 & 39.49 \\
		14920 & 18500 & 23.99 \\
		13823 & 12568 & 9.08 \\
		8566 & 5287 & 38.28 \\
		\hline
	\end{tabular}
	\label{tab:Approx_ratio_full}
\end{table}
\begin{table}[h!]
	\centering
	\caption{Optimal objective function values and approximation ratios with NN surrogate input dimension of 8.}
	\begin{tabular}{|c|c|c|}
		\hline
		\textbf{Obj. Val (\$)} & \textbf{Approx. Val (\$)} & \textbf{Approx. Ratio ( \%)} \\
		\hline
		5228 & 5263 & 0.66 \\
		6410 & 6501 & 1.42 \\
		4568 & 4563 & 0.11 \\
		5411 & 5431 & 0.37 \\
		10166 & 10262 & 0.95 \\
		16335 & 16525 & 1.16 \\
		13443 & 14135 & 5.15 \\
		9174 & 9578 & 4.41 \\
		\hline
	\end{tabular}
	\label{tab:Approx_ratio_reduced}
\end{table}

For our numerical experiments, we select the top five paths and the top three vertiports with the highest AAM demands to impose high-level decisions. The NN surrogate is implemented using a Multi-Layer Perceptron (MLP) with three hidden layers of sizes 256, 128, and 64, respectively. The learning rate is set to 0.01, and the MLP is trained for 500 epochs.


\subsection{Impact of different congestion management policies}
Based on two sets of demand profiles, we generate two different sets of flight schedules to evaluate the performance of the two congestion management strategies described in Section \ref{subsec:aggregation}. 
	Recall the performance of a strategy is evaluated using the congestion metric for each airspace volume, as defined in Section \ref{subsec:HLP}. This metric, represented by the cost function $\Delta_{a, \tau} \coloneqq \max(0, \sum_{f \in F} I_{f, \tau} - Q^{A}_a)$, quantifies the number of AAM flights exceeding the airspace capacity.
As shown in Table~\ref{tab:agg_comparison_1}, the implementation of the type 1 congestion management strategy (naive averaging) reduces congestion across the three test days of the first test set by $33.37\%$, $25.72\%$, and $31.12\%$, respectively, relative to the case when no congestion management strategy is applied. Moreover, the type 2 congestion management strategy (accounting for demand similarity) results in slightly greater reductions in congestion, outperforming the type 1 strategy by an additional $1\%$ to $8\%$. A similar trend is observed in test set 2, where the application of either the type 1 or type 2 strategy significantly decreases congestion compared to the case with no strategy applied. Additionally, the type 2 strategy achieves slightly greater congestion reduction than the type 1 strategy. 

However, we give the following caveat. Even though we observed that type 2 consistently outperforms type 1 day-by-day, this is \emph{not} true at the scale of individual time periods. For example, in time period 7 during the third day of the first test set, we observed that type 1 strategy achieved a congestion level of 8 flights, whereas the type 2 strategy only achieved a congestion level of 10 flights. Nevertheless, even at the scope of individual time periods, both strategies outperformed the case without congestion management strategies. We further caution that future work to test, for example, other demand functions, should be conducted.


	As COPs will define rules---including for conflict management, airspace usage, and demand-capacity balancing---to ensure efficient, fair, and equitable operations across fleet operators and given other airspace users \cite{faa2020uam}, a future congestion management strategy could be encoded as a COP. We argue such strategies must be collaboratively developed by stakeholders, which could include fleet operators, PSUs, and the ANSP such as the FAA. While we demonstrated the utility of two types of congestion management strategies on realistic but synthesized AAM flight schedules, the broader contribution of this paper is providing a way for the FAA to understand potential COPs and how those COPs may be evaluated and evolve to meet operational objectives for specific AAM operations with unique demand profiles, vertiport and routing networks, and the potential for unscheduled pop-up AAM flights. 


\begin{table}[ht]
	\centering
	\caption{Congestion metric comparisons for test set 1. Note that the reduction percentages are shown in parentheses, and that the units are in number of AAM flights above capacity (i.e., lower is better).}
	\label{tab:agg_comparison_1}
	\begin{tabular}{cccc}
		\toprule
		\textbf{Day} & \textbf{No Strat.} & \textbf{Strat. 1 (Red. \%)} & \textbf{Strat. 2 (Red. \%)} \\
		\midrule
		1 & 2158 & 1438 (33.4) & \textbf{1415} (34.4) \\
		2 & 1380 & 1025 (25.7) & \textbf{972} (29.6) \\
		3 & 1070 & 737 (31.1) & \textbf{644} (39.8) \\
		\bottomrule
	\end{tabular}
\end{table}

\begin{table}[ht]
	\centering
	\caption{Congestion metric comparisons for test set 2. Note that the reduction percentages are shown in parentheses, and that the units are in number of AAM flights above capacity (i.e., lower is better).}
	\label{tab:agg_comparison_2}
	\begin{tabular}{cccc}
		\toprule
		\textbf{Day} & \textbf{No Strat.} & \textbf{Strat. 1 (Red. \%)} & \textbf{Strat. 2 (Red. \%)} \\
		\midrule
		1 & 1720 & 1025 (40.4) & \textbf{1013} (41.1) \\
		2 & 2015 & 1391 (31.0) & \textbf{1333} (33.9) \\
		3 & 2743 & 1791 (34.7) & \textbf{1730} (36.9) \\
		\bottomrule
	\end{tabular}
\end{table}

\subsection{Impact of unscheduled, pop-up demand}	\label{ssec:unscheduled_demand}

The analysis of pop-up, unscheduled demand is critical to the operational feasibility of the proposed congestion management concepts. Off-nominal situations, such as prioritized emergency vehicles as well as air traffic not controlled by PSUs (e.g., legacy air traffic such as helicopters and general aviation), can occur. It is crucial to understand the degradation in traffic management performance in these scenarios. For this study, we set the pop-up rate at each OD pair to $10\%$, and the 
	departure times of unscheduled flights follow a continuous uniform distribution,
with the lower and upper bounds corresponding to the start and end of the LLP solve timeframe. Due to the randomness involved in generating the pop-up demand, we generate 100 independent test cases with pop-up demands for each time period.

\begin{figure*}[!htbp]
	\centering
	\includegraphics[width=0.9\linewidth]{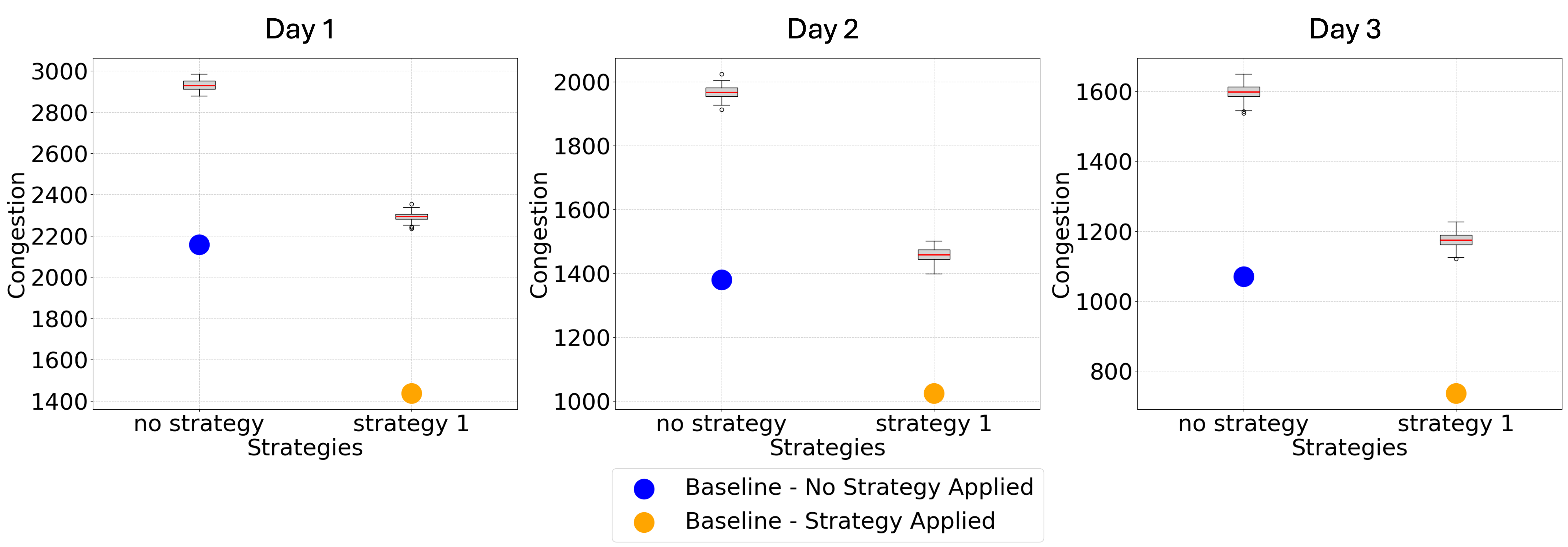}
	\caption{Comparison of resultant congestion with pop-up flight demand. Pop-up flight scenarios shown via box plots, and the baseline scenario represents no pop-up flights (demand profiles from test set 1). Note the $y$-axis is the congestion metric, $\Delta_{a, \tau}$, measured in the number of AAM flights that exceed airspace capacity (i.e., lower is better).}
	\label{fig:ondemand_test1}
\end{figure*}

\begin{figure*}[!htbp]
	\centering
	\includegraphics[width=0.9\linewidth]{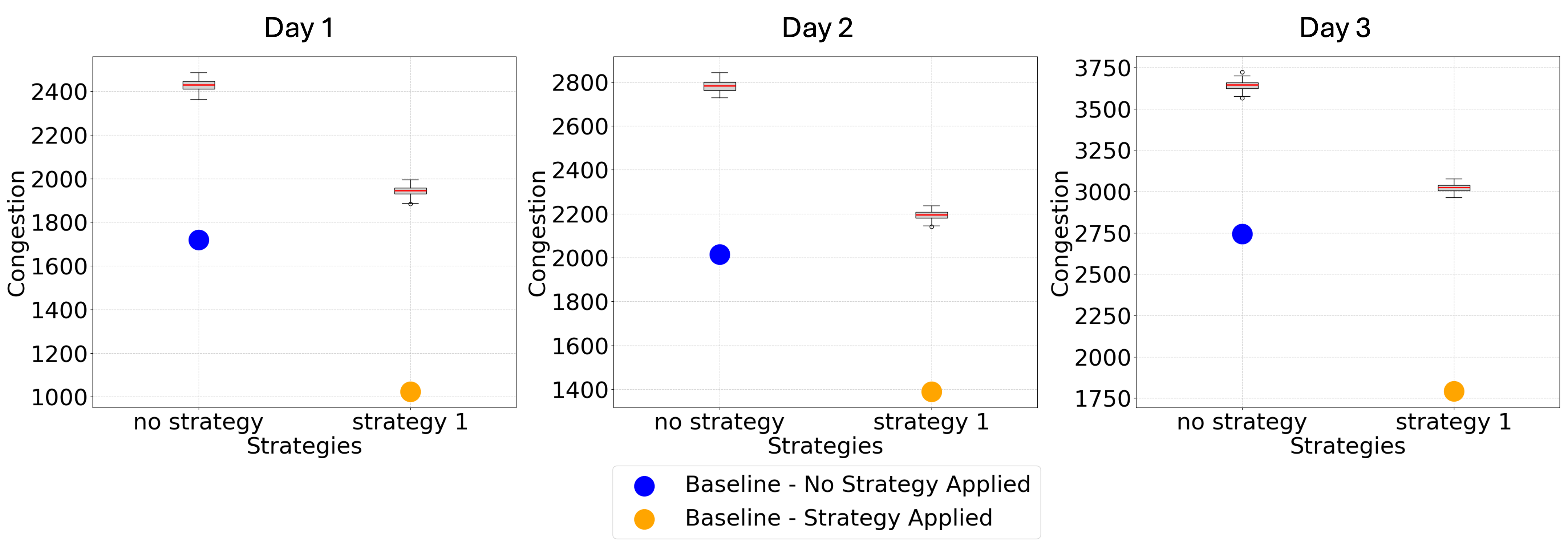}
	\caption{Comparison of resultant congestion with pop-up flight demand. Pop-up flight scenarios shown via box plots, and the baseline scenario represents no pop-up flights (demand profiles from test set 2). Note the $y$-axis is the congestion metric, $\Delta_{a, \tau}$, measured in the number of AAM flights that exceed airspace capacity (i.e., lower is better).}
	\label{fig:ondemand_test2}
\end{figure*}

Figure~\ref{fig:ondemand_test1} illustrates the impact of pop-up demands over three days, 
	with baseline values representing the congestion (as measured by the congestion metric $\Delta_{a, \tau}$) in the absence of pop-up demands.
The results indicate that pop-up demands increase congestion regardless of whether congestion management strategies are implemented. However, when comparing cases with no strategy to those using the type 1 strategy, the overall congestion in scenarios where the congestion management strategies are implemented is significantly lower, even with the increased congestion caused by pop-up demands. Figure~\ref{fig:ondemand_test2} presents similar results, demonstrating that the superior performance of the congestion management strategies holds across different demand profiles and pop-up flight instances. 

\section{Concluding Remarks}	\label{sec:conclusion}

In this work, we contribute an optimization-based approach to evaluate different congestion management strategies within the context of AAM. Our modeling framework mirrors current concepts of operations revolving around xTM: in lieu of the air navigation service provider directly providing strategic traffic flow management capabilities, fleet operators instead interact with federated PSUs in terms of scheduling and routing AAM flights. The role of the air navigation service provider is to monitor congestion levels within the AAM network, and occasionally intervene to achieve specific congestion management objectives. Using a bi-level optimization approach, we are able to capture the two interactions: One between the fleet operator and the PSUs, and another between the PSUs and the air navigation service provider (which in turn impacts the fleet operators). 

We examined the behavior of different congestion management strategies under stochastic, time-varying demand. Using parameterized demand functions, we tested a variety of demand profiles with differing peak periods throughout the day. We also evaluated  demand uncertainties related to pop-up, unscheduled flights. This enabled us to model prioritized AAM flights (e.g., those carrying out emergency response functions)--priority flights must be assigned the shortest paths and without delay. Additionally, these pop-up flights could represent legacy air traffic (e.g., helicopter traffic and general aviation) which are not handled by xTM, but still utilize airspace resources. 


One limitation of work is that we do not know what real-world AAM demand profiles will look like. However, we could have used, e.g., existing rideshare data to better understand demand patterns within a specific metropolitan area. One direction of future work could be to explore a broader range of demand functions, some parametric (such as the one used in this paper) and some non-parametric (such as one derived from real-world demand data). This would help to further validate whether certain congestion management strategies work better than others, and under what conditions. Additionally, we did not consider other impacts on fleet operator actions (e.g., weather conditions). Future studies could  include environmental uncertainties in addition to demand uncertainties.

\section*{NOTICE}

This work was produced for the U.S. Government under Contract 693KA8-22-C-00001 and is subject to Federal Aviation Administration Acquisition Management System Clause 3.5-13, Rights In Data-General (Oct. 2014), Alt. III and Alt. IV (Oct. 2009).

The contents of this document reflect the views of the author and The MITRE Corporation and do not necessarily reflect the views of the Federal Aviation Administration (FAA) or the Department of Transportation (DOT).  Neither the FAA nor the DOT makes any warranty or guarantee, expressed or implied, concerning the content or accuracy of these views.

For further information, please contact The MITRE Corporation, Contracts Management Office, 7515 Colshire Drive, McLean, VA  22102-7539, (703) 983-6000.

\textbf{\copyright 2025 The MITRE Corporation. All Rights Reserved.}

Approved for Public Release; Distribution Unlimited. Case Number 25-0196.


\bibliographystyle{IEEEtran} 
\bibliography{main.bib}

\section*{Author Biographies}

\small{ {\bf Haochen Wu} a Ph.D. Candidate under the supervision of Professor Max Li at the Aerospace Engineering department at University of Michigan, Ann Arbor. He received his master’s degree in Transportation Engineering in 2022 from the University of California, Berkeley, and his B.E. degree in Air Traffic Management in 2020 from the Civil Aviation University of China.
\vspace{1 mm}}

\small{ {\bf Lesley A. Weitz} is a Senior Principal Aerospace Engineer and Chief Scientist in the Transportation Automation Evolution Department in The MITRE Corporation's Center for Integrated Transportation (CIT). Lesley received her PhD in Aerospace Engineering from Texas A\&M University in 2009. 
\vspace{1 mm}}
	
\small{ {\bf Jeffrey M. Henderson} is an Operations Research Analyst in the Transportation Automation Evolution Department in The MITRE Corporation's Center for Integrated Transportation (CIT). Jeffrey received his PhD in Civil and Environmental Engineering from Virginia Tech in 2008.
\vspace{1 mm}}

\small{ {\bf Max Z. Li} is an Assistant Professor of Aerospace Engineering at the University of Michigan, Ann Arbor. He also has courtesy appointments in Civil and Environmental Engineering as well as Industrial and Operations Engineering. Max received his PhD in Aerospace Engineering from the Massachusetts Institute of Technology in 2021. He received his MSE in Systems Engineering and BSE in Electrical Engineering and Mathematics, both from the University of Pennsylvania, in 2018. 
\vspace{1 mm}}


\flushend
\end{document}